\documentclass[11pt]{article}
\usepackage{amssymb}
\usepackage{amsthm}

\usepackage{graphicx}
\usepackage[english]{babel}
\usepackage{lineno}

\usepackage{float}
\usepackage{enumitem}
\usepackage{amsmath}
\usepackage{cleveref}
\usepackage{multirow}
\usepackage{tabularx}
\usepackage{xcolor}
\usepackage{url}
\usepackage{mathtools}

\newtheorem{theorem}{Theorem}[section]

\newcommand{\EQ}{\begin{equation}}
\newcommand{\EN}{\end{equation}}

\newcommand{\Z}{\mathbb{Z}}

\newcommand{\pr}{\indent{\bf Proof. \ }}

\newcommand{\bx}{{\boldsymbol x}}
\newcommand{\bh}{{\boldsymbol h}}

\newcommand{\by}{{\boldsymbol y}}

\newcommand{\bb}{{\boldsymbol b}}
\newcommand{\ba}{{\boldsymbol a}}

\newcommand{\rank}{\mbox{\rm rank}}

\newcommand{\C}{{\cal C}}

\begin{document}




\title{On Sylvester-type constructions of Hadamard matrices
and their modifications}

\author{Dmitrii Zinoviev, Victor Zinoviev}

\maketitle

\begin{abstract}
Using the ideas of concatenation construction of codes over the
$q$-ary alphabet, we modify the known generalized
Sylvester-type construction of the Hadamard matrices. The new construction
is based on two collections of the Hadamard matrices. In particular
this construction involves $m$ Hadamard matrices of order $k$ and $k$ Hadamard matrices
of order $m$. These matrices are not necessary different.
As a result we obtain a Hadamard matrix of order $km$. The new construction
gives many possibilities for construction of the new Hadamard matrices with different
ranks and dimension of kernel.
\end{abstract}



\section{Introduction}

Let $E=\{0,1, \ldots, q-1\}$. A $q$-ary block code $\C$ is an arbitrary nonempty
subset of $E^n$. We say that $\C$ is a $(n,N,d)_q$-code, where $n$ is the
{\em length} of the code, $N$ is the number of codewords, i.e. the
{\em cardinality} of $\C$ and $d$ is the
{\em minimum Hamming distance},
\[
d = d(\C) = \min\{d(\bx, \by):\;\bx \neq \by, \;\bx, \by \in \C\}
\]
where, for $\bx=(x_1, \ldots, x_n)$ and $\by = (y_1, \ldots, y_n)$ from $E^n$,
\[
d(\bx, \by) = |\{j:\, x_j \neq y_j,\;j=1, \ldots, n\}|.
\]

A Hadamard matrix $H$ of order $n$ is an $n \times n$ matrix of
$+1$'s and $-1$'s such that $HH^T=nI$, where $I$ is the $n\times
n$ identity matrix. We know that if a Hadamard matrix $H$ of order
$n$ exists, then $n$ is 1, 2 or a multiple of 4 [1]. Two
Hadamard matrices are equivalent if one can be obtained from
another one by permuting rows and/or columns and multiplying rows and/or
columns by $-1$. We can change the first row and column of $H$
into $+1$'s and we obtain an equivalent Hadamard matrix which is
called normalized. If $+1$'s are replaced by 0's and $-1$'s by
1's, and the complementary of all the rows are added, we obtain a
binary $(n,2n,n/2)$-code which is called a (binary) Hadamard code
and is denoted by $H_n$ [1]. We always assume that the
Hadamard matrix is normalized, so the corresponding Hadamard code
contains the all-zero codeword.

Two structural properties of the binary codes (in particular, the Hadamard codes) are the rank and
the dimension of its kernel. The rank of a binary code $C$ is simply the
dimension of the linear span, $\langle C \rangle$,  of $C$.
The kernel of a binary code $C$ of length $n$ is defined as [2]
\EQ\label{def:kernel}
\ker(C)=\{\textbf{x}\in \Z_2^n : \textbf{x}+C=C \}.
\EN
If the all-zero vector belongs to $C$, then $\ker(C)$ is a linear subcode of $C$.
Note also that if $C$ is linear, then $\ker(C)=C=\langle C \rangle$.
We denote the rank of a binary code $C$ as $\rank(C)$, and the dimension
of the kernel as $\dim(\ker(C))$.
The rank and the dimension of the kernel can be used to distinguish
between nonequivalent binary codes that contain the all-zero vector,
since equivalent ones have the same rank and dimension of the kernel [3].

The goal of the present paper is to describe a new general construction of the binary
Hadamard codes (or matrices). We extend our construction [4], to modify the known generalized
Sylvester-type construction of Hadamard matrices. The new construction is based on two collections of
the Hadamard matrices. In particular it involved $m$ Hadamard matrices of order $k$ and $k$ Hadamard matrices
of order $m$. These matrices are not necessarily different.
As a result we obtain a Hadamard matrix of order $km$. This construction
gives many possibilities for construction of new Hadamard matrices with different
ranks and dimension of kernel.

\section{New constructions}
\label{sec:Construction}

The purpose of this section is to describe a new
construction of the Hadamard matrices, which modifies the known
generalized Sylvester-type construction given by No and Song
[5]. The new construction of the Hadamard matrices
can be considered as a modification of the generalized Sylvester-type
construction from the paper by No and Song [5].

First, recall that by the Sylvester construction of Hadamard matrices from the matrices
$H_n =[h_{i,j}]$ and $H_m$, we mean matrix $H_{mn}$, denoted
by $H_{mn} = H_n \otimes H_m$, obtained
by replacing every element $h_{i,j}$ with the matrix
$h_{i,j} H_m$ (clearly in this case $h_{i,j} \in \{\pm 1\}$),
or with the matrix $h_{i,j} + H_m$ (clearly in this case
$h_{i,j} \in \{0, 1\}$), where, for $H_m = [h_{r,k}]$ and
any fixed pair $i,j$,
\[
h_{i,j} + H_m = [h_{i,j}+h_{r,k}].
\]

In [5], No and Song suggested the following generalized Sylvester-type
construction of Hadamard matrices. We give their result here.

\begin{theorem}\label{general-sylvester} $([5])$
Suppose we have $m$ Hadamard matrices $B_1$, $B_2$, \ldots, $B_m$
of order $k$ (which are not necessary distinct) and a Hadamard
matrix $C=[c_{i,j}]$ of order $m$ where all matrices are over
$\{0,1\}$. Then the matrix $H$,
\[
H=\left[
\begin{array}{cccc}
\;c_{1,1}+B_1\;&\;c_{1,2}+B_2\;&\;\ldots\;&\;c_{1,m}+B_m\;\\
\;c_{2,1}+B_1\;&\;c_{2,2}+B_2\;&\;\ldots\;&\;c_{2,m}+B_m\;\\
       :       &        :      &    :     &      :        \\
\;c_{m,1}+B_1\;&\;c_{m,2}+B_2\;&\;\ldots\;&\;c_{m,m}+B_m\;\\
\end{array}
\right],
\]
is a Hadamard matrix of order $mk$.
\end{theorem}

For a binary vector $\ba$ and $e \in \{0,1\}$ denote for the shortness
the following vector
\[
\ba + e = \ba + e(1,1, \ldots,1).
\]

Our modification of this construction can be stated in the following:

\begin{theorem}\label{modified-general-sylvester}
Suppose we have $m$ Hadamard matrices $A_1$, $A_2$, \ldots, $A_m$
of order $k$ (which are not necessary distinct) and $k$  Hadamard
matrices $B_1$, $B_2$, \ldots, $B_k$ of order $m$ (which are also
not necessary distinct). Let $\ba^{(j)}_i$, $j = 1,2,\ldots,m$,\;
$i = 1,2,\ldots,k$, denote the $i$-th row of $A_j$, let
$B_u = [b^{(u)}_{r,s}]$, $u = 1,2,\ldots,k$,\;
$r,s = 1,2,\ldots,m$ and $\bb^{(u)}_r$ denote $r$-th row of $B_u$.
Then the matrix $H$,
\[
H=\left[
\begin{array}{ccccccc}
&\ba^{(1)}_1+b^{(1)}_{1,1}&\;\ba^{(2)}_1+b^{(1)}_{1,2}& \ldots &\;\ba^{(m)}_1+b^{(1)}_{1,m}&\\
&\ba^{(1)}_2+b^{(2)}_{1,1}&\;\ba^{(2)}_2+b^{(2)}_{1,2}& \ldots &\;\ba^{(m)}_2+b^{(2)}_{1,m}&\\
& \cdots                  &  \cdots                   & \cdots &  \cdots                   &\\
&\ba^{(1)}_k+b^{(k)}_{1,1}&\;\ba^{(2)}_k+b^{(k)}_{1,2}& \ldots &\; \ba^{(m)}_k+b^{(k)}_{1,m}\\\hline

&\ba^{(1)}_1+b^{(1)}_{2,1}&\;\ba^{(2)}_1+b^{(1)}_{2,2}& \ldots\;&\;\ba^{(m)}_1+b^{(1)}_{2,m}&\\
&\ba^{(1)}_2+b^{(2)}_{2,1}&\;\ba^{(2)}_2+b^{(2)}_{2,2}& \ldots\;&\;\ba^{(m)}_2+b^{(2)}_{2,m}&\\
& \cdots                  &  \cdots                   & \cdots      &  \cdots              &\\
&\ba^{(1)}_k+b^{(k)}_{2,1}&\;\ba^{(2)}_k+b^{(k)}_{2,2}&\ldots\;&\;\ba^{(m)}_k+b^{(k)}_{2,m}&\\\hline
&       \cdots       &        \cdots     &   \cdots     &     \cdots                       &\\\hline
&\ba^{(1)}_1+b^{(1)}_{m,1}&\;\ba^{(2)}_1+b^{(1)}_{m,2}& \ldots &\;\ba^{(m)}_1+b^{(1)}_{m,m}&\\
&\ba^{(1)}_2+b^{(2)}_{m,1}&\;\ba^{(2)}_2+b^{(2)}_{m,2}& \ldots &\;\ba^{(m)}_2+b^{(2)}_{m,m}&\\
& \cdots                  &  \cdots                   & \cdots &  \cdots                   &\\
&\ba^{(1)}_k+b^{(k)}_{m,1}&\;\ba^{(2)}_k+b^{(k)}_{m,2}& \ldots &\; \ba^{(m)}_k+b^{(k)}_{m,m}\\

\end{array}
\right],
\]
is a Hadamard matrix of order $mk$.
\end{theorem}

We provide an independent proof of this result in terms of the matrices $A_i$ and $B_j$:

\pr Let $\bh_{i_1}$ and $\bh_{i_2}$ be two different rows of $H$. We have
to consider three different cases.\\
(i) Both rows $\bh_{i_1}$ and $\bh_{i_2}$ belong to the same $i$-th shell of $H$,
which is obtained from the elements $b^{(s)}_{i,r}$, where $s = 1,2, \ldots,  k$  and
$r = 1,2, \ldots, m$. We have for his case
\[
d(\bh_{i_1} \bb_{i_2}) = \sum_{r=1}^m \frac{k}{2} \cdot d(\bb^{(s)},\bb^{(s')}) = \frac{km}{2}.
\]
since $s \neq s'$.\\
(ii) Both rows $\bh_{i_1}$ and $\bb_{i_2}$ belong to the different
$i$-th and $i'$-th shells of $H$, but have the same parameter $s$.
In this case we obtain (taking into account, that $\ba^{(j)}_s = \ba^{(j)}_{s'}$)
\begin{eqnarray*}
d(\bh_{i_1}, \bb_{i_2})
&=& \sum_{r=1}^m d(\ba^{(r)}_s + b^{(s)}_{i,r}, \ba^{(r)}_s + b^{(s)}_{i',r}) \\
&=& k \times \sum_{r=1}^m d(b^{(s)}_{i,r}, b^{(s)}_{i',r})\\
&=& k \times d(\bb^{(s)}_i, \bb^{(s)}_{i'})\\
&=& k \times \frac{m}{2} = \frac{km}{2}.
\end{eqnarray*}
(iii) Both rows $\bh_{i_1}$ and $\bb_{i_2}$ belong to the different
$i$-th and $i'$-th shells of $H$ and correspond to the different
rows $\ba^{(r)}_s$ and $\ba^{(r)}_{s'}$,\, $s\neq s'$.
Now we have
\begin{eqnarray*}
d(\bh_{i_1}, \bb_{i_2})
&=& \sum_{r=1}^m d(\ba^{(r)}_s + b^{(s)}_{i,r}, \ba^{(r)}_{s'} + b^{(s')}_{i',r}) \\
&=& k \times \sum_{r=1}^m d(\ba^{(r)}_s, \ba^{(r)}_{s'})\\
&=& k \times \frac{m}{2}.
\end{eqnarray*}
Indeed, the second equality follows since the following equality
is valid for any different $s$ and $s'$
\[
d(\ba^{(r)}_s, \ba^{(r)}_{s'}) = d(\ba^{(r)}_s, \ba^{(r)}_{s'} + (1,1, \cdots,1)),
\]
implying that
\[
d(\ba^{(r)}_s + b^{(s)}_{i,r}, \ba^{(r)}_{s'} + b^{(s')}_{i',r}) = d(\ba^{(r)}_s, \ba^{(r)}_{s'})
\]
\qed

\section{Acknowledgements}
The research was carried out at the IITP RAS within the framework of fundamental research on
the topic mathematical foundations of the theory of error-correcting codes and at the expense of the
National Science Foundation of Bulgaria (NSFB) under project number 20-51-18002.

\end{document}